\documentclass[preprint,12pt]{elsarticle}
\usepackage{amsmath,amssymb,amsthm}
\usepackage{rotating}
\usepackage{hyperref}
\newtheorem{theorem}{Theorem}

\journal{Discrete Mathematics}

\begin{document}
\begin{frontmatter}

\title{A tight single-change covering design with block size \texorpdfstring{$6$}{6}\tnoteref{ded}}
\tnotetext[ded]{Dedicated to the memory of Donald A.\ Preece (1939--2014),
who posed this problem.}

\author{Richard Bean}
\ead{r.bean1@uq.edu.au}
\address{Cyber Research Centre, University of Queensland, 4072, Australia}

\begin{abstract}
We give a tight single-change covering design with $v=26$ and $k=6$. This
answers Problem~1 of the Nineteenth British Combinatorial Conference, which
asked whether such a design exists with block size greater than $5$. We also
describe the satisfiability search that found the design, including negative
search results at the smallest order not excluded by known bounds, $v=21$.
\end{abstract}

\begin{keyword}
single-change covering design \sep tight single-change covering design \sep
persistent pair \sep satisfiability \sep constraint programming
\MSC[2020] 05B40 \sep 05B30 \sep 68R07
\end{keyword}

\end{frontmatter}

\section{Introduction}
\label{sec:intro}

Following Wallis, Yucas and Zhang~\cite{wyz} and Preece et al.~\cite{tsccd95},
a \emph{tight single-change covering design}, denoted $\mathrm{tsccd}(v,k)$, is
an ordered sequence of $k$-subsets of $S=\{1,2,\dots,v\}$, with $v>k$, such that
\begin{enumerate}
\item[(i)] any two elements of $S$ occur together in at least one block
(covering);
\item[(ii)] each block after the first is obtained from the previous block by
removing one element and inserting another (single change); and
\item[(iii)] the element inserted to form any block after the first has not
previously occurred in a block with any of the other elements of that block
(tightness).
\end{enumerate}
We display a design as an array whose columns are the blocks,
with unchanged elements kept in the same row and printed as dots; the
$\mathrm{tsccd}(26,6)$ constructed here is shown in Table~\ref{tab:design}.
Each insertion of an element, including its first appearance, is a
\emph{transfer}.

Tightness gives a useful counting test. When an element is inserted, it is
paired with the $k-1$ elements already present, none of which it has met before.
Thus the first block creates $\binom{k}{2}$ pairs, each later block creates
exactly $k-1$ new pairs, and no pair is created twice. Since all pairs must be
covered, every pair occurs in exactly one run of consecutive blocks. If the
design has $b$ blocks, then
\begin{equation}
\binom{v}{2}=\binom{k}{2}+(b-1)(k-1),
\qquad\text{so}\qquad
b=\frac{\binom{v}{2}-\binom{k}{2}}{k-1}+1 .
\label{eq:b}
\end{equation}
Counting the $k$ elements of the first block as transfers, the total number of
transfers is $T=b+k-1$. Let $t_i$ be the number of elements transferred exactly
$i$ times, so $\sum_i t_i=v$ and $\sum_i i\,t_i=T$. An element transferred only
once occurs in a single run of $v-k+1$ consecutive blocks; no run can be longer.
We call such an element an \emph{anchor}. For $k\ge3$, a pair occurring together in $v-k$
consecutive blocks is a \emph{persistent pair}~\cite{persistent,tsccd95}; there
are at most $\lfloor k/2\rfloor$ such pairs, and they are disjoint~\cite{tsccd95}.

A $\mathrm{tsccd}(v,k)$ is \emph{standardised}~\cite{tsccd95,p20} if its first
block is $\{1,2,\dots,k\}$, the other elements are first introduced in the order
$k+1,k+2,\dots,v$, and the elements of the first block are first removed in the
order $k,k-1,\dots,1$. Every design can be standardised by relabelling. We use standardised designs
throughout to avoid searching equivalent labellings.

The problem goes back to Nelder's work on the economical formation of
correlation matrices~\cite{nelder} and to the successive-block schemes of Gower
and Preece~\cite{gowerpreece}. Single-change covering designs were introduced
in~\cite{wyz}; tight designs were studied in detail in~\cite{tsccd95}, with a
narrative account in~\cite{insidestory}. The cases $k=2$ (trivial),
$k=3$~\cite{tsccd95}, and $k=4$~\cite{pp12,tsccd95} are well understood. For
$k=5$, Phillips~\cite{p20} found the smallest example, a $\mathrm{tsccd}(20,5)$.
Preece then asked whether any tight single-change covering design exists with
$k>5$; this became Problem~1 of the Nineteenth British Combinatorial
Conference~\cite{bcc19}.

A necessary condition is $v\ge3(k-1)$~\cite{wyz,tsccd95}. For $k=6$,
this gives $v\ge15$, while \eqref{eq:b} requires $v\equiv0$ or $1\pmod5$.
Zhang's bound~\cite{zhang}, $v(v-1)\ge(6v-7k)(k-1)$, becomes
$(v-10)(v-21)\ge0$. Together these conditions leave the orders
\[
v=21,\ 25,\ 26,\ 30,\ 31,\ 35,\ 36,\dots,
\]
the smallest being $v=21$, with $b=40$. We show the following.

\begin{theorem}\label{thm:main}
There exists a tight single-change covering design $\mathrm{tsccd}(26,6)$.
\end{theorem}

\begin{proof}[Proof of Theorem~\ref{thm:main}]
The blocks in \ref{app:blocks} were checked by a separate program.
There are $63$ blocks, each containing six distinct elements of
$\{1,\dots,26\}$, and each consecutive pair of blocks differs by one element.
All $\binom{26}{2}=325$ pairs occur, each in a single run of consecutive
blocks. Thus the sequence satisfies covering, single change and tightness.
The dot array in Table~\ref{tab:design} gives the same sequence.
\end{proof}

The design has $b=63$ blocks, $T=68$ transfers, and
$(t_1,t_2,t_3,t_4)=(6,5,8,7)$. It is shown in Table~\ref{tab:design}, with the
blocks listed in \ref{app:blocks}. 
Section~\ref{sec:model} gives the
satisfiability model, Section~\ref{sec:calib} tests it on the known $k=4$ and
$k=5$ cases, and Section~\ref{sec:k6} reports the $k=6$ search.

\section{A satisfiability formulation}
\label{sec:model}

We represent a standardised $\mathrm{tsccd}(v,k)$ as a system of constraints over
binary variables and use a solver to check if a satisfying assignment exists. The
primary variables are
\[
z_{e,j}\in\{0,1\}\quad(e\in\{1,\dots,v\},\ j\in\{1,\dots,b\}),
\]
where $z_{e,j}=1$ means that element $e$ lies in block $j$. Block sizes are
fixed by $\sum_{e}z_{e,j}=k$ for every $j$. For the single-change condition we introduce
binary variables $s_{e,j}$, with $s_{e,1}=z_{e,1}$ and, for $j\ge2$,
$s_{e,j}=1$ exactly when $e$ is in block $j$ but not in block $j-1$.
We impose $\sum_e s_{e,j}=1$ for each $j\ge2$. Since every block has size $k$,
exactly one element is also removed at each step.

For each unordered pair $\{e,f\}$ and block $j$, we introduce two binary
variables. The first records whether both elements occur in block $j$. The second
marks the beginning of a run of blocks containing the pair. The model requires exactly
one such run for each pair. We also add
redundant constraints that block~$1$ starts $\binom{k}{2}$ pair runs and
each later block starts exactly $k-1$. Standardisation is imposed using
auxiliary variables that record first appearances and first removals.

An element transferred $m$ times occurs in $v-1-(k-2)m$ blocks. At each
transfer it encounters $k-1$ new elements, and each new block of the same
run adds one more.
Between successive runs of an element, the model requires at least $k-2$
blocks where the element is absent. Its first and last appearances must
span at least $v-k+1$ blocks, including both endpoints.
For $e>k$, its first-appearance block $j_e$ satisfies
\[
e-k+1\le j_e\le b-v+k-(v-e).
\]

We generally used the CP-SAT solver of Google OR-Tools~\cite{ortools},
which combines SAT-style clause learning with constraint propagation.
Each unsuccessful run either terminates with an \texttt{INFEASIBLE} status, which rules out a
solution to the encoded constraints, or results in a timeout.
We also tested an equivalent integer linear program with
Gurobi~\cite{gurobi}, and a separate backtracking program that constructs
single-change block sequences directly.

\section{Tests on known designs}
\label{sec:calib}

We first tested the formulation on solved cases. With the first block
fixed by standardisation, CP-SAT found a $\mathrm{tsccd}(12,4)$ ($b=21$) in
about two seconds, and the smaller $k=3$ designs in well under a second. The
backtracking search returned them in milliseconds. However, Gurobi found no
$\mathrm{tsccd}(12,4)$ within $600$ seconds.
For $\mathrm{tsccd}(16,4)$ ($b=39$), adding redundant constraints reduced
the CP-SAT search time from $115$ to $43$ seconds.

For $k=5$, CP-SAT did not find a $\mathrm{tsccd}(20,5)$ within an hour
without fixed anchor runs. In the design of Phillips~\cite{p20} used in our
tests, the five anchors $\{1,6,11,19,20\}$ have length-$16$ runs starting in
blocks $1,2,7,15,16$, respectively. Fixing these runs and asking CP-SAT to
complete the design returned a valid $\mathrm{tsccd}(20,5)$ in $247$ seconds. As a check on the
formulation, fixing all entries of the published design as constraints
returned a satisfying assignment in $0.2$ seconds.
For $k=6$, we similarly fixed selected anchor runs and searched for a
completion.

\section{Block size 6}
\label{sec:k6}

\subsection{The order \texorpdfstring{$v=21$}{v=21}}
\label{sec:v21}

We first searched the smallest order allowed by the bounds above, $\mathrm{tsccd}(21,6)$
($b=40$), using the strategy above. Several hundred CP-SAT runs, on a laptop and
on the University of Queensland Bunya cluster over about a week, did not produce
a design. Many constrained instances ended with solver reports of infeasibility
rather than timeouts.

Call the introduction schedule \emph{front-loaded} if every element first
appears as early as standardisation allows, that is, element $e>k$ first appears
in block $e-k+1$. The $\mathrm{tsccd}(20,5)$ used above has this schedule.
For $v=21$, CP-SAT reported infeasibility after $2421$ seconds with this
schedule fixed and no prescribed transfer counts. To probe nearby schedules,
let the \emph{introduction delay} of a design be
\[
D=\sum_{e>k}\bigl(j_e-(e-k+1)\bigr),
\]
where $j_e$ is the block in which $e$ first appears. Thus $D=0$ precisely for a
front-loaded design. CP-SAT showed that no $\mathrm{tsccd}(21,6)$ exists with
$D\le 10$. In addition to the $D=0$ search, separate runs ruled out
$D\in[0,4]$, $[5,6]$, $[7,8]$ and $[9,10]$. The first three of these intervals
were checked with two random seeds. All the fixed anchor configurations tried
at $v=21$ were also ruled out. Searches for completions of $100$ standardised
openings of $16$ blocks found no design. Searches without fixed introductions
or anchor runs, including searches with prescribed transfer counts, remained
unresolved at their time limits.

These searches leave the existence of a $\mathrm{tsccd}(21,6)$ with $D>10$
unresolved. We next considered $v=25$ and $v=26$.

\subsection{The orders \texorpdfstring{$v=25$ and $v=26$}{v=25 and v=26}}
\label{sec:v25v26}

At $v=25$ ($b=58$) and $v=26$ ($b=63$), we searched with prescribed
anchor runs. A design with $k=6$ can have at most three persistent pairs,
compared with at most two for $k=5$. We tried selected configurations of six
anchors, fixing each anchor $e$ to a run of $v-5$ blocks starting in block~$1$
if $e=1$, and in block $e-5$ otherwise.

Configurations placing three anchors on the last three
elements $\{v-2,v-1,v\}$ were infeasible at each order tried through $v=36$
($8$ seconds at $v=25$, rising to about an hour at $v=36$). The configuration
with anchors $\{1,7,13,19,24,25\}$ at $v=25$, modelled on the
$\mathrm{tsccd}(20,5)$, was ruled out in $264$ seconds. In all, $18$ of the $21$ two-persistent-pair configurations tried at
$v=25$ were refuted, several only after $2$ to $18$ hours of computation.

The configuration that produced a design had three persistent pairs: at $v=26$,
the anchors were $\{1,7,14,15,25,26\}$, with persistent pairs $(1,7)$,
$(14,15)$ and $(25,26)$.
From this configuration CP-SAT returned a satisfying assignment after $4615$
seconds on a task using four CPU cores.

The resulting $\mathrm{tsccd}(26,6)$ is shown in Table~\ref{tab:design}. Its
once-transferred elements are exactly the six anchors $\{1,7,14,15,25,26\}$.
Every element is introduced as early as standardisation permits, so the design
is front-loaded.

\begin{sidewaystable}
\centering
\caption{A standardised $\mathrm{tsccd}(26,6)$ ($v=26$, $k=6$, $b=63$,
$T=68$, $(t_1,t_2,t_3,t_4)=(6,5,8,7)$). Blocks are columns; an unchanged
element is shown as a dot. The design is given in two halves, blocks $1$--$32$
then blocks $33$--$63$; in the second half the leftmost column states block
$33$ in full.}
\label{tab:design}
\footnotesize
\medskip\par\noindent Blocks $1$--$32$:\par\nobreak
\begin{verbatim}
 1  .  .  .  .  .  .  .  .  .  .  .  .  .  .  .  .  .  .  .  .  6  . 13  .  .  .  9  .  .  .  .
 2  .  .  .  .  .  . 13 14  .  .  .  .  .  .  .  .  .  .  .  .  .  .  .  .  .  .  .  . 10  . 11
 3  .  .  .  .  . 12  .  .  .  .  . 18  .  . 21  . 23 24  . 26  .  .  .  .  .  .  .  .  .  .  .
 4  .  .  9 10  .  .  .  . 15  .  .  .  .  .  .  .  .  .  .  .  .  .  .  .  .  .  .  .  . 16  .
 5  .  8  .  . 11  .  .  .  . 16 17  . 19 20  . 22  .  . 25  .  .  .  .  .  .  .  .  .  .  .  .
 6  7  .  .  .  .  .  .  .  .  .  .  .  .  .  .  .  .  .  .  .  .  8  .  3  2  5  .  4  .  .  .
\end{verbatim}
\medskip\par\noindent Blocks $33$--$63$:\par\nobreak
\begin{verbatim}
 9  .  .  .  .  .  .  .  .  .  . 23  .  .  .  .  .  .  .  .  . 17  .  .  .  .  .  .  .  . 11
11 23  . 12  .  . 18  .  .  .  .  .  .  .  .  .  .  .  .  . 19  .  .  .  .  .  .  .  .  .  .
26  .  .  .  .  .  .  .  . 13  .  . 10 11 12  .  .  . 16  .  .  .  .  .  .  .  .  . 10  8  .
17  . 20  . 19 22  . 21  .  .  .  .  .  .  .  .  .  .  .  .  .  . 20  .  .  .  .  .  .  .  .
25  .  .  .  .  .  .  . 24  .  .  .  .  .  .  .  .  .  .  .  .  .  . 22  .  .  .  .  .  .  .
 4  .  .  .  .  .  .  .  .  .  6  .  .  .  .  5  3  8  .  2  .  .  .  . 13  3  6  5  .  .  .
\end{verbatim}
\end{sidewaystable}

\section{Concluding remarks}
\label{sec:remarks}

The existence of $\mathrm{tsccd}(21,6)$ and $\mathrm{tsccd}(25,6)$ remains
unresolved by these searches. The full spectrum of orders for $k=6$, and the
existence question for $k>6$, remain open.

\section*{Acknowledgements}
This work was supported by resources provided by The University of Queensland
Research Computing Centre's Bunya supercomputer~\cite{bunya}, with funding from
The University of Queensland, Brisbane, Australia.

\appendix
\section{The blocks of the \texorpdfstring{$\mathrm{tsccd}(26,6)$}{tsccd(26,6)}}
\label{app:blocks}
The $63$ blocks, as sets, in order:
{\footnotesize
\begin{verbatim}
 1: {1,2,3,4,5,6}      22: {6,7,14,15,25,26}    43: {6,9,13,18,21,24}
 2: {1,2,3,4,5,7}      23: {6,8,14,15,25,26}    44: {6,13,18,21,23,24}
 3: {1,2,3,4,7,8}      24: {8,13,14,15,25,26}   45: {6,10,18,21,23,24}
 4: {1,2,3,7,8,9}      25: {3,13,14,15,25,26}   46: {6,11,18,21,23,24}
 5: {1,2,3,7,8,10}     26: {2,13,14,15,25,26}   47: {6,12,18,21,23,24}
 6: {1,2,3,7,10,11}    27: {5,13,14,15,25,26}   48: {5,12,18,21,23,24}
 7: {1,2,7,10,11,12}   28: {5,9,14,15,25,26}    49: {3,12,18,21,23,24}
 8: {1,7,10,11,12,13}  29: {4,9,14,15,25,26}    50: {8,12,18,21,23,24}
 9: {1,7,10,11,12,14}  30: {4,9,10,15,25,26}    51: {8,16,18,21,23,24}
10: {1,7,11,12,14,15}  31: {4,9,10,16,25,26}    52: {2,16,18,21,23,24}
11: {1,7,12,14,15,16}  32: {4,9,11,16,25,26}    53: {2,16,19,21,23,24}
12: {1,7,12,14,15,17}  33: {4,9,11,17,25,26}    54: {2,16,17,19,21,24}
13: {1,7,14,15,17,18}  34: {4,9,17,23,25,26}    55: {2,16,17,19,20,24}
14: {1,7,14,15,18,19}  35: {4,9,20,23,25,26}    56: {2,16,17,19,20,22}
15: {1,7,14,15,18,20}  36: {4,9,12,20,25,26}    57: {13,16,17,19,20,22}
16: {1,7,14,15,20,21}  37: {4,9,12,19,25,26}    58: {3,16,17,19,20,22}
17: {1,7,14,15,21,22}  38: {4,9,12,22,25,26}    59: {6,16,17,19,20,22}
18: {1,7,14,15,22,23}  39: {4,9,18,22,25,26}    60: {5,16,17,19,20,22}
19: {1,7,14,15,22,24}  40: {4,9,18,21,25,26}    61: {5,10,17,19,20,22}
20: {1,7,14,15,24,25}  41: {4,9,18,21,24,26}    62: {5,8,17,19,20,22}
21: {1,7,14,15,25,26}  42: {4,9,13,18,21,24}    63: {5,8,11,19,20,22}
\end{verbatim}
}


\end{document}